\documentclass[conference]{IEEEtran}

\usepackage{amssymb,amsthm}
\usepackage[cmex10]{amsmath}
\usepackage[pdftex]{graphicx}
\usepackage{mathtools}
\usepackage{cite}

\def\blfootnote{\xdef\@thefnmark{}\@footnotetext}

\newcommand{\C}{\mathcal{C}}
\newcommand{\LT}{\mathcal{L}}
\newcommand{\E}{\mathbb{E}}
\newcommand{\ed}{\,{\buildrel d \over =}\,}
\newcommand{\cd}{\xrightarrow{d}}
\newcommand{\el}{\mathcal{\partial}}
\newcommand{\ninf}{\nu \to \infty}
\newcommand{\limnu}{\lim_{\ninf}}

\newcommand{\bin}[2]{\begin{pmatrix} #1 \\ #2 \end{pmatrix}}
\newcommand{\pr}[1]{{\mathbb{P}}\left\{ {#1} \right\}}

\newtheorem{thm}{Theorem}[section]
\newtheorem*{thma}{Theorem}

\newtheorem{lem}[thm]{Lemma}
\newtheorem{prop}[thm]{Proposition}

\theoremstyle{remark} 
\newtheorem{rem}{Remark}

\newcommand{\eqn}[1]{\begin{equation} #1 \end{equation}}
\newcommand{\eqan}[1]{\begin{align} #1 \end{align}}

\IEEEoverridecommandlockouts

\begin{document}
\belowdisplayskip=5.5pt
\abovedisplayskip=5.5pt 

\title{Mixing Properties of CSMA Networks\\on Partite Graphs}


\author{\IEEEauthorblockN{Alessandro Zocca}
\IEEEauthorblockA{Dept. of Math. \& Comp. Science\\
Eindhoven University of Technology\\
Den Dolech 2, 5612 AZ\\ Eindhoven, The Netherlands\\
a.zocca@tue.nl}
\and
\IEEEauthorblockN{Sem C. Borst}
\IEEEauthorblockA{Dept. of Math. \& Comp. Science\\
Eindhoven University of Technology\\
Den Dolech 2, 5612 AZ\\Eindhoven, The Netherlands\\
s.c.borst@tue.nl}
\and
\IEEEauthorblockN{Johan S. H. van Leeuwaarden}
\IEEEauthorblockA{Dept. of Math. \& Comp. Science\\
Eindhoven University of Technology\\
Den Dolech 2, 5612 AZ\\ Eindhoven, The Netherlands\\
j.s.h.v.leeuwaarden@tue.nl}
\thanks{This research was financially supported by The Netherlands Organization for Scientific Research (NWO) in the framework of the TOP-GO program and by an ERC Starting Grant.}
}

\maketitle

\begin{abstract}
We consider a stylized stochastic model for a wireless CSMA network. Experimental results in prior studies indicate that the model provides remarkably accurate throughput estimates for IEEE 802.11 systems. In particular, the model offers an explanation for the severe spatial unfairness in throughputs observed in such networks with asymmetric interference conditions. Even in symmetric scenarios, however, it may take a long time for the activity process to move between dominant states, giving rise to potential starvation issues. 

In order to gain insight in the transient throughput characteristics and associated starvation effects, we examine in the present paper the behavior of the transition time between dominant activity states. We focus on partite interference graphs, and establish how the magnitude of the transition time scales with the activation rate and the sizes of the various network components. We also prove that in several cases the scaled transition time has an asymptotically exponential distribution as the activation rate grows large, and point out interesting connections with related exponentiality results for rare events and meta-stability phenomena in statistical physics. In addition, we investigate the convergence rate to equilibrium of the activity process in terms of mixing times.
\end{abstract}


\section{Introduction and motivation}
\label{sec1}

We consider a stylized model for a network of nodes sharing a wireless medium according to a CSMA-type protocol. The network is described by an undirected graph $(V, E)$ where the set of vertices~$V$ represents the various nodes of the network and the set of edges $E \subseteq V \times V$ indicates which pairs of nodes interfere. In other words, nodes that are neighbors in the interference graph are prevented from simultaneous activity, and thus the independent sets correspond to the feasible joint activity states. A node is said to be blocked whenever the node itself or any of its neighbors is active, and unblocked otherwise. Each node activates (starts a transmission) at an exponential rate~$\nu$ whenever it is unblocked. The transmission durations of nodes are independent and exponentially distributed with unit mean.

Let $\Omega^* \subseteq \{0, 1\}^V$ be the collection of incidence vectors of the independent sets of the interference graph $(V, E)$, and let $X_t^* \in \Omega^*$ denote the joint activity state at time~$t$, with its $i$-th element indicating whether or not node~$i$ is active at time~$t$. Then $(X_t^*)_{t \geq 0}$ is a reversible Markov process~\cite{Kelly79} with stationary distribution
\eqn{
\label{statt}
\pi_x(\nu) = \lim_{t \to\infty} \pr{X_t^*=x} =
\frac{\nu^{\|x\|_1}}{\sum_{y \in \Omega^*} \nu^{\|y\|_1}},
\quad x \in \Omega^*,
}
with $\|x\|_1 = \sum_{i \in V} x_i$ the number of active nodes in state~$x$.

The above model has a long history \cite{BKMS87,Kelly85,KBC87}. It was rediscovered in the context of IEEE 802.11 systems in~\cite{WK05}, and further developed in that setting in \cite{DT06,SGMCK08}. While the modeling of the IEEE 802.11 back-off mechanism is less detailed than in the model of Bianchi~\cite{Bianchi00}, the general interference graph offers greater versatility and covers a broad range of topologies. Experimental results in~\cite{LKLW10} demonstrate that the model, while idealized, provides remarkably accurate throughput estimates for actual IEEE 802.11 systems. It is also worth observing that the model amounts to a special instance of a loss network \cite{SR04,ZZ99}, and that the stationary distribution corresponds to the Gibbs measure of the hard-core model in statistical physics~\cite{GS08,Haggstrom97}.

An activity state in $\Omega^*$ is called {\it dominant} if it corresponds to an independent set of maximum size $\max_{x \in \Omega^*} \|x\|_1$. It follows from the stationary distribution~\eqref{statt} that only the dominant states retain probability mass as $\ninf$. This causes a severe problem, called {\it spatial unfairness}, when some nodes belong to fewer maximum-size independent sets than others as a result of asymmetric interference conditions~\cite{VLDJ09}, and receive far lower or even zero throughputs.

Even in symmetric scenarios where the long-term throughputs are equal,
however, potential starvation issues can occur, since it may take
a long time for the activity process to move between dominant states.
Consider for example an interference graph where the nodes can be
partitioned into two independent sets of maximum size.
The activity process will spend roughly half of the time in each of
the two associated dominant states as $\ninf$.
Since each of the nodes is active in one of the dominant states,
they will all receive equal long-term throughputs, so spatial
unfairness is not an issue in the long run.
Yet, another source of severe unfairness arises as $\ninf$, because
it will take an extremely long time for the activity process to move
from one dominant state to the other (resembling {\it meta-stability
phenomena} in statistical physics \cite{DenHollander04,OS95}).
As a result, each node will experience long sequences of transmissions
in rapid succession, interspersed with extended periods of starvation.

In order to gain insight in the above issue, we explore in the
present paper the behavior of the Markov process $(X_t^*)_{t \geq 0}$
in a scenario with a complete $K$-partite interference graph $(V,E)$.
By this we mean that the nodes can be partitioned into $K$~disjoint
sets called {\it components}, such that two nodes block each other,
if and only if they belong to different components.
The key contributions of the present paper may be summarized as follows:

(i) We examine the asymptotic behavior (as $\ninf$) of the transition
time between various activity states, exploiting classical results for
absorption times in birth-and-death processes \cite{KMcG59,K65,K71}
and a representation of the transition time as a geometric random sum.

(ii) We establish that the magnitude of the transition time scales
as $\Theta(\nu^{\beta-1})$, reflecting how `rigid' the dominant
activity pattern is as function of the activation rate~$\nu$
and the exponent~$\beta$, which is completely determined by the
sizes of the various components.

(iii) We prove that in several cases the scaled transition time has an asymptotically
exponential distribution and compare it with related exponentiality results for rare events \cite{A82,AB92,AB93,GR05,K66,K79}.

(iv) We investigate the rate of convergence of the activity process to the equilibrium distribution, and demonstrate that the mixing time is of the same order as the escape time of the second-largest component.

The remainder of the paper is organized as follows.
In Section~\ref{sec2} we present a model description and some
preliminary results.
In Section~\ref{sec3} we examine the asymptotic behavior of the
transition time in the case of a bipartite graph.
Section~\ref{sec4} describes the extension of the results to
arbitrary partite graphs.
In Section~\ref{sec5} we investigate the rate of convergence to
equilibrium in terms of mixing times.
Section~\ref{sec6} concludes with some remarks and a review of
further extensions.

\section{Model description and preliminary results}
\label{sec2}

Consider the Markov process $(X_t^*)_{t \geq 0}$ as described in the
introduction with a complete $K$-partite interference graph $(V,E)$.
Thus the nodes in~$V$ can be partitioned into $K$~disjoint sets
called {\it components}, such that two nodes are connected by
an edge in~$E$, if and only if they belong to different components.
In view of the symmetry, all the states with the same number of
active nodes in a given component can be aggregated, and we only
need to keep track of the number of active nodes, if any,
and the component they belong to.
This \textit{state aggregation} yields an equivalent Markov process
$(X_t)_{t \geq 0}$ on a star-shaped state space~$\Omega$ with
$K$~branches which emanate from a common root node and correspond
to the components of the interference graph.
Specifically, $\Omega =
\{0\} \cup \{(k,l) : 1 \leq l \leq L_k, \, 1 \leq k \leq K\}$,
where the center state~0 indicates that none of the nodes is active
and state $(k,l)$ corresponds to the situation where $l$~nodes are
active in the $k$-th component, denoted by $\C_k$, with $L_k$ denoting the size of $\C_k$ (see Figure 1). Notice that $\Omega$ can be alternative described as $\{0\} \cup \bigcup_{k=1}^K \C_k$. The transition rates of the process $(X_t)_{t \geq 0}$ are given
by $q(0, (k, 1)) = L_k \nu$, $q((k,l), (k,l+1)) = (L_k - l) \nu$,
$l = 1, \dots, L_k - 1$, $q((k,1), 0) = 1$, and $q((k,l), (k,l-1)) = l$, $l = 2, \dots, L_k$, $k = 1, \dots, K$.
\vspace{-0.25cm}
\begin{figure}[hbtp]
\begin{center}
\includegraphics[height=3.3cm]{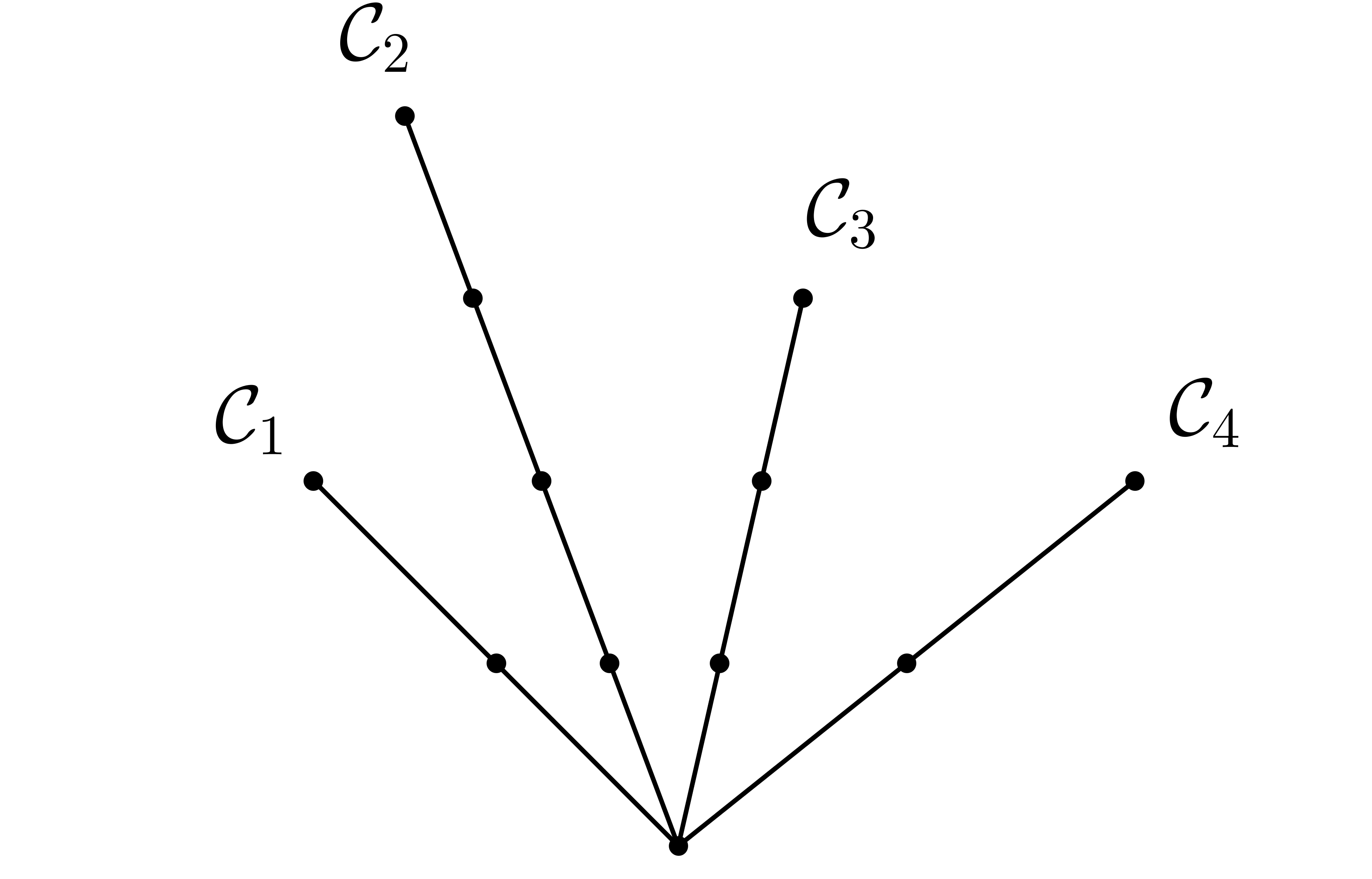}
\vspace{-0.3cm}
\caption{Example of the state space $\Omega$ with $K=4$.}
\vspace{-0.5cm}
\end{center}
\end{figure}

The stationary distribution of the process $(X_t)_{t \geq 0}$ reads
\eqan{
\label{statdist}
\pi_0(\nu) &= \Big(1 + \sum_{k = 1}^{K} \sum_{l = 1}^{L_k}
\bin{L_k}{l} \nu^l\Big)^{-1}, \\
\pi_{(k,l)}(\nu) &= \pi_0(\nu) \bin{L_k}{l} \nu^l, \quad l = 1, \dots, L_k, \, k = 1, \dots, K. \nonumber
}

Denote by
\[
T_{(k_1,l_1), (k_2,l_2)}(\nu)=
\inf\{t > 0: X_t = (k_2,l_2) | X_0 = (k_1,l_1)\}
\]
the first-passage time of state $(k_2,l_2)$ starting in $(k_1,l_1)$, which we also refer to as the {\it transition time} from state $(k_1,l_1)$ to state $(k_2,l_2)$.
In the next sections we will analyze the asymptotic behavior of the
transition time as $\ninf$.
As discussed in the introduction, the transition time provides
useful insight in transient throughput characteristics
and starvation effects in wireless CSMA networks.

In preparation for the asymptotic analysis of general scenarios,
we first present a few results for the case where the two states
$(k_1,l_1)$ and $(k_2,l_2)$ belong to the same component, i.e.
$k_1 = k_2$ and $l_1 > l_2$.
The presence of the other components is not relevant then for the
transition time, and hence we focus on the case of just a single
component ($K = 1$), and drop the component index for now.
When $K = 1$, the process $(X_t)_{t \geq 0}$ evolves as an elementary
birth-and-death process on the state space $\{L, L-1, \dots, 1,0\}$,
so that we can exploit several classical results for such processes.

\subsection{Asymptotic growth rate}

We first state how the expectation of the transition time
$T_{l_1,l_2}(\nu)$ scales as $\ninf$.
Here and throughout the paper we write $f(\nu) \sim g(\nu)$ to
indicate that $\limnu f(\nu) / g(\nu) = 1$ as $\ninf$ for any two
real-valued functions $f(\cdot)$ and $g(\cdot)$.
\vspace{-0.05cm}
\begin{prop}
\label{mean}
For $L \geq l_1 > l_2 \geq 0$,
\[
\E T_{l_1,l_2}(\nu) \sim \frac{l_2!(L-l_2-1)!}{L!} \nu^{L-l_2-1},
\quad \ninf.
\]
\end{prop}
\begin{IEEEproof}
First observe that
$\E T_{l_1,l_2}(\nu) = \sum_{l=l_1}^{l_2+1} \E T_{l,l-1}(\nu)$,
so we can exploit a general result for birth-and-death
processes~\cite{K65}, which in the present case says that,
for $0 < l \leq L$,
$
\E T_{l,l-1}(\nu) =
\frac{1}{l} \sum_{n=l}^{L} \frac{\pi_n(\nu)}{\pi_l(\nu)}.
$
Now~\eqref{statdist} implies that
$\pi_n(\nu) = {\rm o}(\pi_L(\nu))$ as $\ninf$ for all
$n = l, \dots, L - 1$, so that
\[
\E T_{l,l-1}(\nu) \sim \frac{1}{l} \frac{\pi_{L}(\nu)}{\pi_{l}(\nu)} =
\frac{(L-l)!(l-1)!}{L!} \nu^{L-l}, \quad \ninf.
\]
Thus $\E T_{l,l-1}(\nu) = {\rm o}(\E T_{l_2+1,l_2})$ as
$\ninf$ for all $l = l_1,\dots, $ $l_2$, and hence $\E T_{l_1,l_2}(\nu) \sim \E T_{l_2+1,l_2}(\nu)$ as $\ninf$.
\end{IEEEproof}

In order to gain insight in starvation effects, we are particularly
interested in the time for the activity process to reach the center
state~0, referred to as {\it escape time}, because at such points
in time nodes in other components have an opportunity to activate.
Proposition~\ref{mean} shows that
\eqan{
\label{sdp}
\E T_{l_1,0}(\nu) \sim \frac{1}{L} \nu^{L-1}, \quad \ninf.
}
Hence, the escape time grows as a power of~$\nu$, where the
exponent corresponds to the component size minus one, and is
asymptotically not influenced by the starting state~$l_1$.

\subsection{Asymptotic exponentiality}

Having obtained the asymptotic growth rate, we now turn attention
to the distribution of the scaled escape time, and show that it
has an asymptotically exponential distribution.
We will leverage the following well-known result for birth-and-death
processes, which is commonly attributed to Keilson~\cite{K71} or
Karlin and McGregor~\cite{KMcG59}.

\begin{thm}
\label{thm:km}
Consider a birth-and-death process with generator matrix~$Q$ on the state space $\{0,\dots,d\}$ started at state~$d$. Assume that $0$ is an absorbing state, and that the other birth rates $\{\lambda_i\}_{i=1}^{d-1}$ and death rates $\{\mu_i\}_{i=1}^d$ are positive. Then the absorption time in state~$0$ is distributed as the sum of $d$~independent exponential random variables whose rate parameters are the $d$~nonzero eigenvalues of $-Q$.
\end{thm}

Let $Q(\nu)$ be the generator matrix of the birth-and-death process $(X_t)_{t \geq 0}$ on the state space $\{L,L-1,\dots,1,0\}$, with 0 an absorbing state. Let
$0 < \alpha_1(\nu) <  \alpha_2(\nu) < \dots < \alpha_L(\nu) $
denote the non-zero eigenvalues of $-Q(\nu)$. The property 
that the eigenvalues $\{\alpha_i(\nu)\}_{i=1}^L$ are distinct, real and strictly positive, is well-known~\cite{LR54}, see also the proof of Lemma~\ref{asymptotic1}.

Applying Theorem~\ref{thm:km}, we obtain that
\vspace{-0.22cm}\\
\eqn{
\label{sumexp}
T_{L,0}(\nu) \ed \sum_{i=1}^{L} Y_i(\nu),
}
\vspace{-0.25cm}\\
where $Y_1(\nu),\dots,Y_L(\nu)$ are independent exponentially distributed random variables with $\E Y_i(\nu)=1/\alpha_i(\nu)$. 

Intuitively, the eigenvalue $\alpha_1(\nu)$ will become really small as $\ninf$, and so the mean escape time $\E T_{L,0}(\nu)$ will be dominated by $\E Y_1(\nu)= 1/\alpha_1(\nu)$. Combining~\eqref{sdp} and~\eqref{sumexp} suggests that
\vspace{-0.22cm}\\
\[
\alpha_1(\nu) \sim \frac{1}{\E T_{L,0}(\nu)} \sim \frac{L}{\nu^{L-1}}, \quad \ninf.
\]
Indeed, we can make this precise with the following result for the growth rates of the eigenvalues as $\ninf$ and their relation to the mean escape time $\E T_{L,0}(\nu)$. The proof of the result is presented in Appendix~A, and exploits detailed information about the growth rates of the eigenvalues obtained via symmetrization and the Gershgorin circle theorem.

\begin{lem}
\label{asymptotic1}
\[
\limnu \alpha_i(\nu) \cdot \E T_{L,0}(\nu) =
\left\{\begin{array}{ll}
1, & \hbox{$i=1$,} \\
\infty, & \hbox{$i=2,\dots, L$.}
\end{array}\right.
\]
\end{lem}

The above lemma shows that the smallest eigenvalue $\alpha_1(\nu)$
becomes dominant as $\ninf$, but also proves the asymptotic exponentiality of the escape time. Indeed \eqref{sumexp} means that the Laplace transform of the scaled escape time is
\vspace{-0.05cm}
\eqn{
\LT_{T_{L,0}(\nu)/ \E T_{L,0}(\nu)}(s)=\prod_{i=1}^L \Big(1+\frac{s}{\alpha_i(\nu) \cdot \E T_{L,0}(\nu)} \Big)^{-1}.
 \nonumber }
Lemma~\ref{asymptotic1} then implies that  
\vspace{-0.05cm}
\[
\limnu \LT_{T_{L,0}(\nu) / \E T_{L,0}(\nu)}(s) = \frac{1}{1+s}.
\]
The continuity theorem for Laplace transforms then yields that the scaled escape time has an asymptotically exponential distribution as stated in the next theorem, where ${\rm Exp}(\lambda)$ denotes an exponentially distributed random variable with mean $1/\lambda$.
\vspace{-0.1cm}
\begin{thm}
\label{thm:expo}
$ \displaystyle
\frac{T_{L,0}(\nu)}{\E T_{L,0}(\nu)} \cd {\rm Exp}(1), \quad \ninf.
$
\end{thm}

The above result may be interpreted as follows.
For large~$\nu$, state~$L$ is frequently visited,
while state~0 becomes rare.
This suggests that the probability of hitting state~0 before the
first return to state~$L$ becomes small.
So the time $T_{L,0}(\nu)$ consists of a geometric number of
excursions from~$L$ which return to~$L$ without hitting~0,
followed by part of the excursion that hits~0.
Hence, apart from this final excursion, $T_{L,0}(\nu)$ is the sum
of a large geometric number of i.i.d.~random variables,
which should be asymptotically exponential.

The fact that the time until the first occurrence of a rare event
is asymptotically exponential, is a widely observed phenomenon in
probability~\cite{K79}.
In order to establish exponentiality of hitting some subset~$B$ of
the state space, one typically decomposes the process into
regenerative cycles, and assumes that (i) the probability of
hitting~$B$ in a single cycle is small, and (ii) the length of the
cycle in which $B$ is hit is asymptotically negligible compared
with the mean cycle length \cite{GR05,K79}.
For the case $K=1$, both assumptions hold, and an alternative proof
of Theorem~\ref{thm:expo} can be obtained using~\cite[Thm.1]{GR05}
(which is a generalized version of the theorem proved in~\cite{K66}).
However, this general theorem for regenerative processes will fail
for the majority of cases considered in this paper.
For instance, in the case $K=2$ and $L_1=L_2$ (two identical
components), the process $(X_t)_{t \geq 0}$ will exhibit
{\it bistable behavior}, in the sense that, as $\ninf$, it spends
extremely long periods circling around either one of the leaf nodes.
For analyzing the transition time between the two dominant leaf nodes,
assumption (i) remains valid, because it is still extremely hard to
go from one leaf to the other, but assumption (ii) will be violated.
Indeed, when the process crosses all the way from one leaf to the other,
the length of the latter cycle will be of the same order as the
entire hitting time.
Despite the fact that the regenerative approach no longer works,
we will prove (using different methods) in the next section that
the scaled transition time between dominant states is still
asymptotically exponential.

Let us finally remark that for reversible Markov chains similar
exponentiality results were established in \cite{A82,AB92,AB93}.
Aldous~\cite{A82} showed that a result like Theorem~\ref{thm:expo}
can be expected when the underlying Markov process converges
rapidly to stationarity.
Indeed, the Markov process $(X_t)_{t \geq 0}$ for the case $K=1$
turns out to have a small mixing time.
However, the picture changes drastically in the case $K \geq 2$,
for which we prove in Section~\ref{sec5} that the mixing time is
$\Theta(\nu^{\beta-1})$ with $\beta$ the size of the second-largest
component;
hence the mixing time is of the same order as the escape time of
the second-longest branch, see~\eqref{sdp}.
For large~$\nu$, this implies that the Markov process is extremely
slowly mixing, which is another way of understanding why the
regenerative approach sketched above fails in this situation.

\section{Bipartite interference graph}
\label{sec3}

In this section we examine the asymptotic behavior of the transition time $T_{(k_1,l_1), (k_2,l_2)}(\nu)$ as $\ninf$ for any pair of states $(k_1,l_1)$ and $(k_2,l_2)$ when the interference graph is bipartite.
We will establish how the expectation scales (Theorem~\ref{thm:meanbipartite}) and use two different methods to prove that the scaled transition time has an asymptotically exponential distribution when states $(k_1,l_1)$ and $(k_2,l_2)$ belong to different branches (Theorem~\ref{thm:bipartiteprob}).
Most of these results in fact extend to arbitrary partite graphs, as will be shown in the next section.
However, we treat the case of a bipartite graph separately first, since it allows us to develop the key ideas in a relatively transparent setting, while being sufficiently rich to exhibit the essential qualitative characteristics of general scenarios.

Consider a bipartite interference graph, i.e.~$K=2$, where the two components $\C_1$ and $\C_2$ have sizes $L_1$ and $L_2$, respectively.
In this case, the two branches of the state space along with the root node form a simple linear array.
For notational convenience, we will therefore relabel the states, indicating the number of active nodes in $\C_1$ by a negative integer in $\{-L_1,\dots,-1,\}$ and representing the number of those in $\C_2$ by a positive integer in $\{1,\dots,L_2\}$.

The activity process $(X_t)_{t \geq 0}$ on $\Omega=\C_1 \cup \{0\} \cup \C_2$ evolves as a birth-and-death process with transition rates
\eqan{
q(l,l+1)&=
\begin{cases}
|l| & \text{ if } -L_1 \leq l<0,\\
(L_2-l) \nu & \text{ if }0 \leq l < L_2,
\end{cases}
\nonumber \\
q(l,l-1)&=
\begin{cases}
(L_1-|l|) \nu & \text{ if } -L_1 < l \leq 0, \\
l & \text{ if } 0 < l \leq L_2,
\end{cases}
\nonumber
}
yielding
\eqan{
\pi_{- l}(\nu) &= \pi_0(\nu) \bin{L_1}{l} \nu^l,
\quad l = 1, \dots, L_1, \label{statdistleft} \\
\pi_{m}(\nu) &= \pi_0(\nu) \bin{L_2}{m} \nu^m,
\quad m = 1, \dots, L_2. \label{statdistright}
}

\subsection{Asymptotic growth rate}

We first determine how the expectation of the transition time
$T_{-l_1,l_2}(\nu)$ scales as $\ninf$.
\begin{thm}
\label{thm:meanbipartite}
For any $0 < l_1 \leq L_1$ and $0 < l_2 \leq L_2$,
\eqn{ \label{meanl1l2}
\E T_{-l_1,l_2}(\nu) \sim  \frac{L_1 + L_2}{L_1 L_2} \nu^{L_1-1},
\quad \ninf.
}
\end{thm}

In the CSMA context, the above theorem implies that when at least one
of the nodes in component~$k$ is currently active, $k = 1, 2$, it will
approximately take an amount of time
$\frac{L_1 + L_2}{L_1 L_2} \nu^{L_k-1}$ before any number of nodes
in the other component $3 - k$ will have a chance to transmit.

\begin{IEEEproof}
As before, note that for any $l_1 < l_2$ $\E T_{l_1,l_2}(\nu) = \sum_{l=l_1}^{l_2-1} \E T_{l,l+1}(\nu)$, so we can again use a general result for birth-and-death processes~\cite{K65}, which in the present case says that
$
\E T_{l,l+1}(\nu) =
\frac{1}{q(l,l+1)} \sum_{n=-L_1}^{l} \frac{\pi_{n}(\nu)}{\pi_{l}(\nu)}
$, for $-L_1 \leq l < L_2$.\\
Now~\eqref{statdistleft} and~\eqref{statdistright} imply
that $\pi_n(\nu) = {\rm o}(\pi_{-L_1}(\nu))$ as $\ninf$ for all
$n = - L_1 + 1, \dots, \min\{L_1 - 1, L_2\}$
and $\pi_n(\nu) = {\rm o}(\pi_l(\nu))$ as $\ninf$ for all
$n = - L_1, \dots, l - 1$, when $l = L_1 + 1, \dots, L_2$.
Thus, for $-L_1 \leq l < 0$, as $\ninf$,
\[
\E T_{l,l+1}(\nu) \sim \frac{1}{l} \frac{\pi_{-L_1}(\nu)}{\pi_{l}(\nu)} =
\frac{(L_1-|l|)!(|l|-1)!}{L_1!} \nu^{L_1-|l|},
\]
and, for $0 \leq l \leq L_2$, as $\ninf$, $\E T_{l,l+1}(\nu)$ scales as
\eqn{\label{meanlpositive}
\begin{cases}
\frac{\nu^{-1}}{(L_2-l)} \frac{\pi_{-L_1}(\nu)}{\pi_{l}(\nu)} =  
\frac{(L_2-l-1)! \, l!}{L_2!} \nu^{L_1-l-1} & \text{if } l < L_1, \\
\frac{\nu^{-1}}{(L_2-L_1)} \frac{\pi_{-L_1}(\nu) + \pi_{L_1}(\nu)}{\pi_{L_1}(\nu)} =  
\frac{L_2!+L_1! (L_2-L_1)!}{(L_2-L_1)L_2!}\nu^{-1} & \text{if } l = L_1, \\
\frac{\nu^{-1}}{(L_2-l)} \frac{\pi_l(\nu)}{\pi_l(\nu)} = 
\frac{1}{L_2-l} \nu^{-1} & \text{if } l > L_1.
\end{cases}
}
It follows that $\E T_{l,l+1}(\nu) = {\rm o}(\E T_{- 1, 0}(\nu))$
as $\ninf$ for all $l = - L_1, \dots, - 2$ and $l = 1, \dots, L_2 - 1$.
Hence, as $\ninf$,
\eqn{
\E T_{-l_1,l_2}(\nu) \sim \E T_{-1,0}(\nu) + \E T_{0,1}(\nu) \sim
\frac{1}{L_1} \nu^{L_1-1}+ \frac{1}{L_2} \nu^{L_1-1}, \nonumber
}
which yields~\eqref{meanl1l2}.
\end{IEEEproof}

\begin{rem} \label{rem1}
Theorem~\ref{thm:meanbipartite} only deals with the case where the
two states $-l_1$ and $l_2$ belong to different components, but extends
to the case where the two states belong to the same component.
In case $0 \leq l_2 < l_1 \leq L_1$, the transition time
$T_{-l_1,-l_2}(\nu)$ is not influenced by the presence of the second
component, and its expected value thus follows from
Proposition~\ref{mean}.
In case $0 \leq l_1 < l_2 \leq L_2$, it may be deduced
from~\eqref{meanlpositive} that $\E T_{l_1,l_2}(\nu) \sim \E T_{l_1,l_1+1}(\nu)$ as $\ninf$, yielding
\[
\E T_{l_1,l_2}(\nu) \sim \begin{cases}
\frac{(L_2-l_1-1)! \, l_1!}{L_2!} \nu^{L_1-l_1-1} & \text{if } l_1 < L_1 \\
\frac{L_2! +L_1! (L_2-L_1)!}{(L_2-L_1) L_2!}\nu^{-1} & \text{if } l_1 = L_1 \\
\frac{1}{L_2-l_1} \nu^{-1} & \text{if } l_1 > L_1
\end{cases}, \: \ninf.
\]
\end{rem}

\subsection{Asymptotic exponentiality}

In the previous subsection we obtained the asymptotic growth rate
of the mean transition time $\E T_{-l_1,l_2}(\nu)$ as $\ninf$.
We now turn attention to the scaled transition time and will prove
that it has an asymptotically exponential distribution as stated in
the next theorem.

\begin{thm}
\label{thm:bipartiteprob}
For any $0 < l_1 \leq L_1$ and $0 < l_2 \leq L_2$
\[
\frac{T_{-l_1,l_2}(\nu)}{\E T_{-l_1,l_2}(\nu)} \cd {\rm Exp}(1),
\quad \ninf.
\]
\end{thm}

Before providing a proof, we first present an interpretation of the above theorem.
In order for the process to make a transition from state $-l_1$ to state $l_2$, it must first reach state~0. Proposition~\ref{mean} indicates that the expected transition time from state $-L_1$ to state $-l_1$ is asymptotically negligible compared to the expected transition time from state $-l_1$ to state~0. This means that the transition time from state $-l_1$ to state~0 asymptotically behaves as the transition time from state $-L_1$ to state~0. Theorem~\ref{thm:expo} shows that the latter time, after scaling, has an asymptotically exponential distribution. 
Once the process has reached state~0 for the first time, there are two possible scenarios. With probability $p = L_1 / (L_1 + L_2)$, the process returns to state $-1$, and then almost surely falls back to state $-L_1$ quite rapidly (compared with the time to reach state~0). With probability $1 - p$, the process moves to state~1, and then most likely is attracted to state $l_2$ quite quickly.

In conclusion, in order for the process to make a transition from state $-l_1$ to state $l_2$ it must visit state~0 once plus a geometrically distributed number of times with parameter~$p$. The successive time periods to reach state~0 are independent and, after scaling, exponentially distributed. Observing that the sum of one plus a geometrically distributed number of independent and exponentially distributed random variables is again exponentially distributed, we arrive at the statement of the above theorem.

Obviously, the above argumentation is heuristic, and in the next two subsections we will provide a rigorous proof. The first proof method is analytical in nature, and relies on an asymptotic characterization of the eigenvalues of the generator matrix, while the second proof is more probabilistic, and in fact closely mirrors the above intuitive explanation.

\begin{rem} \label{rem2}
Theorem~\ref{thm:bipartiteprob} only deals with the case where the
two states $-l_1$ and $l_2$ belong to different components, but
extends to the case where the two states belong to the same component.
In case $0 \leq l_2 < l_1 \leq L_1$, the transition time
$T_{-l_1,-l_2}(\nu)$ is not affected by the presence of the second
component, and thus has an asymptotically exponential distribution
according to Theorem~\ref{thm:expo}.
In contrast, in case $0 \leq l_1 < l_2 \leq L_2$, a fundamentally
different situation emerges.
In that case, a series of upward transitions from state~$l_1$ to
state~$l_2$ occurs in rapid succession with high probability,
so that the transition time converges in distribution to~0,
even though its expectation may not tend to~0 and in fact grow to
infinity as $\ninf$ when $l_1 < L_1$, as observed in Remark~\ref{rem1}.
\end{rem}

\subsection{Analytical approach}

We now extend the analytical approach of Section \ref{sec2} to
prove the asymptotic exponentiality of the scaled transition time
as stated in Theorem~\ref{thm:bipartiteprob}.
For convenience, we restrict attention to $l_1 = L_1$,
but the proof readily extends to any $l_1 \in \{1, \dots, L_1 - 1\}$. 

Let $Q(\nu)$ be the generator matrix of the birth-and-death process $(X_t)_{t \geq 0}$ on the state space $\{-L_1, \dots, -1, 0, 1, \dots, l_2\}$ with $l_2$ an absorbing state. Let $0 < \gamma_1(\nu) <  \gamma_2(\nu) < \dots < \gamma_{L_1+l_2}(\nu)$ denote the non-zero eigenvalues of $-Q(\nu)$.

Theorem~\ref{thm:km} implies that $T_{-L_1,l_2}(\nu) \ed \sum_{i=1}^{L_1+l_2} Y_i(\nu)$, where $Y_1(\nu),\dots,Y_{L_1+l_2}(\nu)$ are independent and exponentially distributed random variables with $\E Y_i(\nu) = 1/\gamma_i(\nu)$.

The following lemma shows that the smallest eigenvalue $\gamma_1(\nu)$ becomes dominant as $\ninf$, and implies the asymptotic exponentiality stated in Theorem~\ref{thm:bipartiteprob}. The proof is similar to that of Lemma~\ref{asymptotic1} and thus omitted.

\begin{lem}
\label{asymptotic2}
\eqan{
\limnu \gamma_i(\nu) \cdot \E T_{-L_1,l_2}(\nu) =
\left\{\begin{array}{ll}
1, & \hbox{$i=1$,} \\
\infty, & \hbox{$i=2,\dots,L_1+l_2.$}
\end{array} \right.
\nonumber}
\end{lem}

\subsection{Probabilistic approach}

We now present an alternative, probabilistic approach to establish the asymptotic exponentiality of the scaled transition time as stated in Theorem~\ref{thm:bipartiteprob}. In contrast to the analytical method as used in the previous subsection, the probabilistic approach in fact extends to arbitrary partite graphs, but we first focus on the case of a bipartite graph in order to illuminate the key ideas. The approach relies on a stochastic decomposition of the transition time into independent random variables which are easier to handle. In order to obtain the stochastic decomposition, we consider the evolution of the process as it makes a transition from state $- l_1$ to a state $l_2$, and define the following random variables:

$\bullet$ $T_{-l_1,-1}(\nu)$:
time to reach state $-1$ for the first time;

$\bullet$ $T_{-1,0}^{(0)}(\nu)$:
time to reach state~0 once the process has reached state $-1$ for the first time;

$\bullet$ $N$:
number of times the process makes a transition $0 \to -1$ before the first transition $0 \to 1$ occurs;

$\bullet$ $\hat{T}_{0,-1}^{(i)}(\nu)$:
time spent in state~0 before the $i$-th transition back to state $-1$, $i = 1, \dots, N$;

$\bullet$ $\hat{T}_{0,1}(\nu)$:
time spent in state~0 before the first transition to state~1;

$\bullet$ $T_{-1,0}^{(i)}(\nu)$:
time to return to state~0 after the $i$-th transition back to state $-1$, $i = 1, \dots, N$;

$\bullet$ $T_{1,l_2}(\nu)$:
time to reach state~$l_2$ once the process has reached state~1 for the first time.

By definition, the transition time may be represented as
\eqan{
\label{birepr}
T_{-l_1,l_2} \ed & T_{-l_1,-1} + T_{-1,0}^{(0)} \\
&+ \sum_{i=1}^{N} \left(\hat{T}_{0,-1}^{(i)} + T_{-1,0}^{(i)}\right) +
\hat{T}_{0,1} + T_{1,l_2}, \nonumber
}
where the dependence on the parameter~$\nu$ is suppressed for compactness. We can make the following observations:\\
$\bullet$
The random variables $T_{-l_1,-1}$ and $T_{1,l_2}$ are distributed as typical hitting times for the respective pairs of states;\\
$\bullet$
The random variables $T_{-1,0}^{(i)}$, $i = 0, 1, \dots, N$, are i.i.d.\ copies of a typical transition time $T_{-1,0}$;\\
$\bullet$
The random variables $\hat{T}_{0,1}$ and $\hat{T}_{0,-1}^{(i)}$, $i = 1, 2, \dots, N$, are i.i.d.\ copies of a random variable $T_0 \ed {\rm Exp}((L_1+L_2) \nu)$, which represents the residence time in state~0;\\
$\bullet$
$N \ed {\rm Geo}(p)$, $p=\frac{L_1}{L_1+L_2}$, independent of the parameter~$\nu$;\\
$\bullet$
All the random variables are independent.

Based on the stochastic representation in~\eqref{birepr} and the above observations, we now proceed to give a proof of the asymptotic exponentiality of the transition time as stated in Theorem~\ref{thm:bipartiteprob}.
The proof consists of three main parts: 
(i) the first part establishes that the distribution of the scaled
transition time $T_{-l_1,l_2}(\nu)$ asymptotically coincides with
that of a dominant term $U(\nu)$, which involves a random sum of
i.i.d.\ random variables (Proposition~\ref{asymbounds});
(ii) the second part serves to identify the asymptotic
behavior of each of these random variables (Proposition~\ref{ael0});
(iii) the third part then shows that a scaled random sum
asymptotically behaves as the random sum of the scaled terms
(Proposition~\ref{asyfinitemean}).

Part~(i):
According to Theorem~\ref{thm:meanbipartite} and Remark~\ref{rem1},
the mean values of the random variables $T_{-1,0}^{(i)}$,
$i = 0, 1, \dots, N$, asymptotically dominate, i.e. 
they are an order-of-magnitude larger than those of all other
random variables in~\eqref{birepr} as $\ninf$.
This suggests that the asymptotic behavior of the transition time
$T_{-l_1,l_2}(\nu)$ will be determined by that of
$
U(\nu) = \sum_{i=0}^{N} T_{-1,0}^{(i)}(\nu),
$
and in particular that the distribution of the scaled transition
time $T_{-l_1,l_2}(\nu) / \E T_{-l_1,l_2}(\nu)$ asymptotically
coincides, as $\ninf$, with that of $U(\nu) / \E U(\nu)$.
This follows as a special case of Proposition~\ref{asymbounds} below,
whose proof is technical, but relatively straightforward,
and omitted because of page limitations.
In passing we note that the above observations also imply
that the expectation of $T_{-l_1,l_2}(\nu)$ scales as
$
\E U(\nu) = \big(1 + \E N(\nu)\big) \, \E T_{-1, 0}(\nu) =\frac{L_1 + L_2}{L_2} \, \E T_{-1, 0}(\nu),
$
with $\E T_{-1, 0}(\nu) \sim \nu^{L_1 - 1} / L_1$ as $\ninf$
according to Proposition~\ref{mean}, which corroborates
Theorem~\ref{thm:meanbipartite}.

Part~(ii):
In order to determine the asymptotic behavior of the random
variable $T_{-1, 0}(\nu)$, we first observe that
$T_{-L_1,0}(\nu) \ed T_{-L_1,-1}(\nu) + T_{-1,0}(\nu).$
According to Theorem~\ref{thm:meanbipartite}, the mean value of
$T_{-L_1,-1}(\nu)$ is asymptotically negligible compared to that of
$T_{-1,0}(\nu)$, which suggests that the asymptotic behavior of
$T_{-1,0}(\nu)$ is equivalent to that of $T_{-L_1,0}(\nu)$.
Since Theorem~\ref{thm:expo} states that
$T_{-L_1,0}(\nu) / \E T_{-L_1,0}(\nu)$ has an asymptotically
exponential distribution, this would imply that the same holds for
the scaled hitting time $T_{-1,0}(\nu) / \E T_{-1,0}(\nu)$, as is covered as a special case of Proposition~\ref{ael0}, stated and proved below.

Part~(iii):
It remains to be shown that $U(\nu) / \E U(\nu)$,
with $U(\nu) = \sum_{i=0}^{N} T_{-1,0}^{(i)}(\nu)$, asymptotically
behaves as the random variable  $\frac{1}{1 + \E N} \sum_{i = 0}^{N} Y_i,$
when $T_{-1,0}(\nu) / \E T_{-1,0}(\nu) \cd Y$, and $Y_1$, $Y_2, \dots$ are i.i.d.\ copies of the random variable~$Y$.
This follows as a special case of Proposition~\ref{asyfinitemean}
below, whose proof is in Appendix~B.
In particular, if $N$ is geometrically distributed and $Y$ is
exponentially distributed, then $U(\nu) / \E U(\nu)$ has an
asymptotically exponential distribution as well.

\begin{prop}
\label{asymbounds}
Let $T(\nu), U(\nu),$ $V(\nu),W(\nu)$ be non-\\negative random variables such that \\
{\rm (i)} $\limnu \E V(\nu)/ \E U(\nu) = 0=\limnu \E W(\nu)/ \E U(\nu)$;\\
{\rm (ii)} For every $\nu > 0$,
$U - V \leq_{\rm st} T \leq_{\rm st} U + W$, i.e.~$\forall \, t > 0$ $\pr{U - V > t} \leq \pr{T > t} \leq \pr{U + W > t}$; \\
{\rm (iii)} $U(\nu) / \E U(\nu) \cd X$ as $\ninf$, where $X$ is a random variable independent of~$\nu$ with continuous c.d.f.. \\
Then $T(\nu) / \E T(\nu) \cd X$ as $\ninf$.
\end{prop}

\begin{prop}
\label{ael0}
For any $ 0 < l \leq L$,
\[
\frac{T_{l,0}(\nu)}{\E T_{l,0}(\nu)} \cd {\rm Exp}(1), \quad \ninf.
\]
\end{prop}
\begin{IEEEproof}
The birth-and-death structure of the process and the strong Markov
property yield the stochastic identity $T_{L,0}(\nu) \ed T_{L,l}(\nu)+T_{l,0}(\nu)$, which gives the
stochastic bounds
$
T_{L,0}(\nu) - T_{L,l}(\nu) \leq_{{\rm st}} T_{l,0}(\nu) \leq_{{\rm st}} T_{L,0}(\nu)
$
(the two terms in the lower bound being dependent).
It follows from Theorem~\ref{thm:expo} that
$T_{L,0}(\nu)/ \E T_{L,0}(\nu) \cd {\rm Exp}(1)$ as $\ninf$.
In order to complete the proof, we can then use Proposition~\ref{asymbounds}, taking $U(\nu)=T_{L,0}(\nu)$,
$V(\nu)=T_{L,l}(\nu)$ and $W(\nu)=0$.
The condition which needs to be checked is
$\limnu \frac{\E V(\nu)}{\E U(\nu)} = 0$, which follows directly
from Proposition~\ref{mean}.
\end{IEEEproof}
\begin{prop}
\label{asyfinitemean}
For every $\nu >0$, define $S_M(\nu):=\sum_{i=1}^{M} X_i(\nu)$, where $M$ is an integer-valued random variable and $\{X_i(\nu)\}_{i \geq 1}$ is a sequence, independent of~$M$, of i.i.d.\ copies of a random variable $X(\nu)$, with $\E X(\nu) < \infty$. Assume that $X(\nu) / \E X(\nu) \cd Y$ as $\ninf$, where $Y$ is some unit-mean random variable. Then
\[
\frac{S_M(\nu)}{\E S_M(\nu)} \cd \frac{1}{\E M} \sum_{i=1}^{M} Y_i,
\quad \ninf,
\]
where $\{Y_i\}_{i \geq 1}$ is a sequence, independent of~$M$, of i.i.d.\ copies of the random variable~$Y$.
\end{prop}

\section{Arbitrary partite graphs}
\label{sec4}

In this section we investigate the asymptotic behavior of the transition time $T_{(k_1,l_1), (k_2,l_2)}(\nu)$ as $\ninf$ for any pair of states $(k_1,l_1)$ and $(k_2,l_2)$ in arbitrary partite graphs. While the key ideas are similar to those used for the bipartite graph in the previous section, there are now some dependencies that require specific treatment and some different qualitative features that arise in certain scenarios.
In particular, it turns out that the scaled transition time may no longer have an asymptotically exponential distribution. This scenario arises when there are longer branches than $k_1$, in which case the number of returns to the root node that dominate the transition time is no longer geometrically distributed plus one, but just geometrically distributed.

As in the case of the bipartite graph, the proof approach involves a stochastic representation of the transition time, but some of the terms are now no longer entirely independent. In order to obtain the stochastic representation, we consider the evolution of the process as it makes a transition from a state $(k_1,l_1)$ to a state $(k_2,l_2)$, and define the following random variables:\\
$\bullet$ $T_{(k_1, l_1), (k_1, 1)}(\nu)$:
time to reach state $(k_1, 1)$ for the first time;\\
$\bullet$ $T_{(k_1, 1), 0}^{(0)}(\nu)$:
time to reach state~$0$ after state $(k_1, 1)$ is visited for the first time;\\
$\bullet$ $N_k$:
number of times the process makes a transition $0 \to (k, 1)$, $k \neq k_2$, before the first transition $0 \to (k_2, 1)$ occurs;\\
$\bullet$ $\hat{T}_{0, (k, 1)}^{(i)}(\nu)$:
time spent in state~0 before the $i$-th transition back to state $(k, 1)$, $k \neq k_2$, $i = 1, \dots, N_k$;\\
$\bullet$ $\hat{T}_{0, (k_2, 1)}(\nu)$:
time spent in state~0 before the first transition to state $(k_2, 1)$;\\
$\bullet$ $T_{(k, 1), 0}^{(i)}(\nu)$:
time to return to state~0 after the $i$-th transition back to state $(k, 1)$, $k \neq k_2$, $i = 1, \dots, N_k$;\\
$\bullet$ $T_{(k_2, 1), (k_2, l_2)}(\nu)$:
time to reach state $(k_2, l_2)$ once the process has reached state $(k_2, 1)$ for the first time.

By definition, the transition time may be represented as
\eqan{ \label{kpartiterepr}
& T_{(k_1, l_1), (k_2, l_2)} \ed
T_{(k_1, l_1), (k_1, 1)} + T_{(k_1, 1), 0}^{(0)} \\
+ & \sum_{k \neq k_2} \sum_{i = 1}^{N_k}
\left(\hat{T}_{0, (k, 1)}^{(i)} + T_{(k, 1), 0}^{(i)}\right)
+ \hat{T}_{0, (k_2, 1)} + T_{(k_2, 1), (k_2, l_2)}, \nonumber
}
where the dependence on the parameter~$\nu$ is suppressed for compactness. 
Define $L := \sum_{k=1}^{K} L_k$, $\gamma := L \nu$, and $p_k := L_k / L$ for all $k = 1, \dots, K$.

We can make the following observations:\\
$\bullet$
The random variables $T_{(k_1,l_1),(k_1,1)}$ and $T_{(k_2,1),(k_2,l_2)}$ are distributed as typical transition times for the respective pairs of states;\\
$\bullet$
The random variables $T_{(k,1), 0}^{(i)}$ are i.i.d.\ copies of a typical transition time $T_{(k,1), 0}$, $i = 0, \dots, N_k$, $k \neq k_2$;\\
$\bullet$
The random variables $\hat{T}_{0, (k_2, 1)}$ and $\hat{T}_{0, (k, 1)}^{(i)}$, $k \neq k_2$, $i = 1, \dots, N_k$, are i.i.d.\ copies of a random variable $T_0 \ed {\rm Exp}(\gamma)$, which is the residence time in state~0;\\
$\bullet$
$N := \sum_{k \neq k_2} N_k \ed {\rm Geo}(1-p_{k_2})$, independent of the parameter~$\nu$;\\
$\bullet$
Given $N = n$, $\bar{N} = (N_1, \dots, N_{k_2-1}, N_{k_2+1}, \dots, N_K)$ has a multinomial distribution with parameters~$n$ and $\bar{p}_1, \dots,$  $\bar{p}_{k_2-1}, \bar{p}_{k_2+1}, \dots, \bar{p}_{K}$, with $\bar{p}_k = p_k / (1 - p_{k_2})$
\eqan{
& \pr{\bar N = (n_1, \dots, n_{k_2-1},n_{k_2+1},\dots, n_{K})}\nonumber \\
& =p_{k_2}(1-p_{k_2})^{\sum_{k \neq k_2} n_k} \bin{\sum_{k \neq k_2} n_k}{n_1,\dots, 
n_K} \prod_{k\neq k_2} \bar{p}_{k}^{n_k} \nonumber \\
&= p_{k_2} \bin{\sum_{k \neq k_2} n_k }{n_1, \dots,
n_K} \prod_{k\neq k_2} p_{k}^{n_k}; \label{multinomial}
}
\vspace{-0.25cm}\\
$\bullet$
All the random variables representing time durations are mutually independent, as well as independent of the random variables $N_k$, $k \neq k_2$. 

Define $L_*:=\max_{j \neq k_2} L_j$, $K_*:=\{k \neq k_2: L_k = L_*\}$,
and $p_*:=|K_*| L_*/ (|K_*|L_* + L_{k_2})$.

We first use the stochastic representation~\eqref{kpartiterepr} to establish how the expectation of the transition time scales.

\begin{thm}
\label{thm:kparmean}
For $0 < l_1 \leq L_{k_1}$, $0 < l_2 \leq L_{k_2}$, $k_1 \neq k_2$,
\[
\E T_{(k_1, l_1), (k_2, l_2)}(\nu) \sim
\left(\frac{I_{\{k_1 \in K_*\}}}{L_*} + \frac{|K_*|}{L_{k_2}}\right)
\nu^{L_* - 1}, \quad \ninf.
\]
\end{thm}

In the CSMA context, the above theorem implies that for any
component~$k$ when some other component is presently active,
it will roughly take an amount of time of the order $\nu^{L^*}$,
$L^* = \max_{j \neq k} L_j$ before any number of nodes in component~$k$
will get an opportunity to transmit. The proof of the above theorem is subsumed in that of the next one. 

We now proceed to determine the asymptotic distribution of the scaled transition time.

\begin{thm}
\label{thm:asybehavior}
For $0 < l_1 \leq L_{k_1}$, $0 < l_2 \leq L_{k_2}$, $k_1 \neq k_2$,
\[
\frac{T_{(k_1, l_1), (k_2, l_2)}(\nu)}{\E T_{(k_1, l_1), (k_2, l_2)}(\nu)} \cd
\frac{1}{\E M} \sum_{i=1}^{M} Y_i, \quad \ninf,
\]
where $M \ed {\rm Geo}(p_*) + I_{\{k_1 \in K_* \}}$ and $Y_i$ are independent exponentially distributed random variables with unit mean.
\end{thm}
\begin{IEEEproof}
The proof is similar to that for the bipartite graph and is based on
the stochastic representation~\eqref{kpartiterepr} of the
transition time $T_{(k_1, l_1), (k_2, l_2)}(\nu)$.
Proposition~\ref{mean} implies that
$\E T_{(k_1,l_1),(k_1,1)}(\nu) \sim \nu^{L_{k_1}-2}$ 
and $\E T_{(k, 1),0}(\nu) \sim \nu^{L_k-1}$ for every branch~$k$.
Moreover, $\E \hat{T}_{0,(k, 1)}(\nu) = {\rm o}(1)$, and it is easily
verified that $\E T_{(k_2,1),(k_2,l_2)}(\nu) = {\rm o}(\nu^{L_* - 1})$
similarly as in the proof of Theorem~\ref{thm:meanbipartite}.
Thus the asymptotically dominant term in~\eqref{kpartiterepr} is
\vspace{-0.2cm}\\
\[
T_{(k_1, 1), 0} (\nu) I_{\{k_1 \in K_* \}} +
\sum_{k \in K_*} \sum_{i = 1}^{N_k} T_{(k, 1), 0}^{(i)}(\nu) \ed
\sum_{i = 1}^{M} T_{(k, 1), 0}^{(i)}(\nu),
\]
\vspace{-0.2cm}\\
for $M \ed {\rm Geo}(p_*) + I_{\{k_1 \in K_*\}}$ and for any
$k \in K_*$, where the latter stochastic equality follows
from~\eqref{multinomial} and the fact that the random variables
$T_{(k,1), 0}$ are identically distributed for all $k \in K^*$.
Taking $X_i (\nu) := T_{(k, 1), 0}^{(i)}(\nu)$ and applying
Propositions~\ref{asymbounds}-\ref{asyfinitemean}
then completes the proof.
Also, noting that $\E M = I_{\{k_1 \in K_*\}} + p^* / (1 - p^*)$
and $\E T_{(k, 1), 0} = \nu^{L^* - 1} / L^*$ yields the statement
of Theorem~\ref{thm:kparmean}.
\end{IEEEproof}

The $k$-th branch, $k \neq k_2$, is called \textit{weakly dominant}
if $k \in K_*$, i.e.~if $L_k \geq L_j$ for all $j\neq k_2$. 
Based on Theorem~\ref{thm:asybehavior}, we may distinguish two scenarios, depending on whether branch $k_1$ is weakly dominant or not.

Suppose that branch $k_1$ is weakly dominant, i.e.~$L_{k_1} = L_*$. In this case $M \ed {\rm Geo}(p_*)+1$ and so Theorem~\ref{thm:asybehavior} implies that for every $0 < l_1 \leq L_{k_1}$ and $0 < l_2 \leq L_{k_2}$
\[
\frac{T_{(k_1, l_1), (k_2, l_2)}(\nu)}{\E T_{(k_1, l_1), (k_2, l_2)}(\nu)} \cd
{\rm Exp}(1), \quad \ninf.
\]
Thus in this case the scaled transition time converges to an exponential random variable with unit mean as $\ninf$.

Suppose instead that branch $k_1$ is not weakly dominant, i.e. $L_{k_1} < L_*$. In this case $M \ed {\rm Geo}(p_*)$ and so Theorem~\ref{thm:asybehavior} implies that for every $0 < l_1 \leq L_{k_1}$ and $0 < l_2 \leq L_{k_2}$
\[
\frac{T_{(k_1, l_1), (k_2, l_2)}(\nu)}{\E T_{(k_1, l_1), (k_2, l_2)}(\nu)} \cd
\frac{1-p_*}{p_*}  \sum_{i=1}^{{\rm Geo}(p_*)} Y_i, \quad \ninf,
\]
where $Y_i$ are independent and exponentially distributed random variables with unit mean.

\section{Mixing times}
\label{sec5}

In the previous sections we have analyzed the transient behavior of
our Markov process $(X_t)_{t \geq 0}$ in terms of hitting times.
In this section we turn attention to the long-run behavior of
our Markov process and in particular examine the rate of
convergence to the stationary distribution.
We measure the rate of convergence in terms of the total variation
distance and explore the intimate connection between the hitting
times and the so-called {\it mixing time} of our Markov process.
The mixing time describes the time required for the distance to
stationarity to become small.
It turns out that the mixing time is largely determined by the time
it takes the process to escape from the second-longest branch in
the state space, as formalized in Theorem~\ref{thm:mix} below.
In this section, we assume without loss of generality that the
branches are indexed such that $L_1 \geq L_2 \geq \dots \geq L_K$,
and denote by $\el_k$ the leaf node of the $k$-th branch.
Also, we attach the starting state $X(0) = x \in \Omega$ as
a superscript to our notation for the Markov process, and thus
write $(X_t^x)_{t\geq 0}$.

\subsection{Main result}

Our objective is to bound the maximal distance over $x\in\Omega$, measured in terms of total variation, between the distribution at time~$t$ and the stationary distribution:
\eqn{
d(t):=\max_{x\in\Omega}\|\pr{X_t^x\in \cdot}-\pi\|_{\rm TV}.\nonumber
}
We define the mixing time of our process as
\eqn{
t_{{\rm mix}}(\epsilon,\nu)=\inf\{t \geq 0 : d(t)\leq \epsilon\}.\nonumber
}

\begin{thm}
\label{thm:mix}
The mixing time of the Markov process $(X_t)_{t\geq 0}$ satisfies
\eqan{
t_{{\rm mix}}(\epsilon, \nu) = \Theta( \nu^{L_2-1}), \nonumber
}
i.e.~$\exists \, C_1(\epsilon), C_2(\epsilon) > 0 \, \, \exists \, \nu_0> 0 $
such that for all $\nu > \nu_0$
\eqan{
C_1(\epsilon) \, \nu^{L_2-1} \leq t_{{\rm mix}}(\epsilon,\nu) \leq C_2(\epsilon) \, \nu^{L_2-1}.
\nonumber}
\end{thm}

Theorem \ref{thm:mix} shows that it can take an extremely long time for the process $(X_t)_{t\geq 0}$ to reach stationarity, especially when $\nu$ is large. Such a long mixing time is typically due to the process being stuck for a considerably period in one of the components, and thus not visiting the states in the other components. This is particularly relevant when in the network has two or more dominant components which together attract the entire probability mass in the limit as $\ninf$. Indeed, in this case the mixing time provides an indication how long it takes for a certain fairness among the dominant components to occur. We will prove Theorem~\ref{thm:mix} by deriving an upper bound for
$d(t)$ using coupling in Subsection~\ref{subsec:up} and a matching
lower bound for $d(t)$ using the bottleneck ratio and the notion of
conductance in Subsection~\ref{subsec:low}.

\subsection{Upper bound using coupling}
\label{subsec:up}

It can be easily established that
\eqn{
d(t)\leq \bar d(t):=\max_{x,y\in\Omega}\|\pr{X_t^x\in \cdot}-\pr{X_t^y \in \cdot}\|_{\rm TV}.
\nonumber}

We consider $\bar d(t)$ instead of $d(t)$, because it can be bounded using a standard coupling technique. Consider all couplings of the processes $(X_t^x,X_t^y)$ with the property that both $(X_t^x)$ and $(X_t^y)$ are Markov processes that have the same generator. Moreover, assume that the coupling is such that the two processes stay together at all times after they have met for the first time.
Under these assumptions, we have that
\eqan{
&\|\pr{X_t^x\in \cdot}-\pr{X_t^y \in \cdot}\|_{\rm TV}\nonumber \\
&=\max_{A\subseteq \Omega}|\pr{X_t^x\in A}-\pr{X_t^y \in A}|\nonumber\\
&=\max_{A\subseteq \Omega}\pr{X_t^x \in A, X_t^y \not\in A}\nonumber \\
&\leq \pr{X_t^x\neq X_t^y}=\pr{\tau_{\rm couple}^{x,y}>t},
\nonumber}
where the {\it coupling time}
$
\tau_{\rm couple}^{x,y}=\min\{t\geq 0:X_t^x=X_t^y\}
$
denotes the first time the two processes $(X_t^x)$ and $(X_t^y)$ meet. Therefore,
$
d(t)\leq \bar d(t)\leq \max_{x,y\in\Omega}\pr{\tau_{\rm couple}^{x,y}>t},
$
and the strength of this coupling inequality depends of course heavily on the choice of the coupling.

We now introduce a birth-and-death process $(M_t)_{t\geq 0}$ that will play a crucial role in our coupling. Let $(M_t)_{t\geq 0}$ describe the position of a particle that lives only on the two longest branches of the state space $\C_1\cup \{0\} \cup \C_2$, starts in the leaf node $\el_2$, and moves within the two branches according to the same transition rates as $(X_t)_{t\geq 0}$. Call the particle that moves according to $(M_t)_{t\geq 0}$ particle~0. Consider also a particle~1 and particle~2 whose positions are governed by the coupled Markov process $(X_t^x,X_t^y)$. To be more specific, we denote the exact position of particle~$i$ at time $t$ by $(K_t^i, L_t^i)$, $i = 0, 1, 2$, with $K_t^i$ the branch and $L_t^i$ the level.

\begin{prop}
\label{prop:couple}
For all $x,y\in\Omega$, the coupling time $\tau_{\rm couple}^{x,y}$, with $x = (k_1, l_1)$
and $y = (k_2, l_2)$ is stochastically bounded from above by the absorption time $T_{\el_2,\el_1}$ of state~$\el_1$ of the continuous-time birth-and-death process $(M_t)_{t\geq 0}$ starting from state $M_0=\el_2$.
\end{prop}

\begin{IEEEproof}
The coupling is such that particles~1 and~2 stay together at all times after they have met for the first time. Before that time, whenever a particle~$i$ ($i=1,2$) resides in $\C_1 \cup \C_2$, it is coupled to particle~0 in the following way: whenever particle~0 moves towards the root, particle~$i$ moves towards the root. In order to construct this coupling, we introduce a Poisson clock with rate $L_1$. When the clock ticks, we first generate a $[0,1]$ uniform random variable~$U$ and then do the following: move particle~$i$ down one level (towards the root) if $U < L_t^i /L_1$; also move particle~0 down one level if $U < L_t^0 /L_1$. This coupling ensures that when $L_t^0 = L_t^i$ ($i=1,2$) and particle~0 moves down one level, then so does particle~$i$.
A consequence of the above coupling is that whenever particle~0 enters $\C_1$ both particle~1 and particle~2 reside in $\C_1$. Extend the coupling by assuming that when particle~0 meets a particle~$i$ ($i=1,2$) in some state in $\C_1$, the two particles keep making the same transitions as long as they are in $\C_1$. Hence, inevitably, by the time particle~0 reaches $\el_1$, the particles~1 and~2 are coupled.
\end{IEEEproof}

Since, for every $\nu >0$ and for all $x,y\in \Omega$,
$
\tau_{\rm couple}^{x,y} (\nu) \leq_{{\rm st}} T_{\el_2,\el_1} (\nu),\nonumber
$
we arrive at the following result:
\eqn{\label{dbound}
d(t)\leq \bar d(t)\leq \pr{T_{\el_2,\el_1}(\nu) >t}\leq \frac{\E T_{\el_2,\el_1} (\nu)}{t}.
}
It then follows immediately that 
$
t_{{\rm mix}}(\epsilon,\nu)\leq \epsilon^{-1}\E T_{\el_2,\el_1}(\nu).
$

Using Theorem \ref{thm:kparmean}, we thus obtain the following upper bounds on the distance to stationarity and the mixing time.
\begin{prop} For the Markov process $(X_t)_{t\geq 0}$ the maximal total variation distance is bounded by
\eqn{
d(t)\leq \frac{\E T_{\el_2,\el_1}(\nu)}{t}\sim \frac{1}{t} \left( \frac{L_1+L_2}{L_1 L_2} \right) \nu^{L_2-1}, \quad \ninf,
\nonumber}
and the mixing time is bounded by
\eqn{
t_{{\rm mix}}(\epsilon,\nu) \leq \frac{1}{\epsilon} \, \E T_{\el_2,\el_1}(\nu) \sim \frac{1}{\epsilon} \left( \frac{L_1+L_2}{L_1 L_2} \right) \nu^{L_2-1}, \, \ninf.
\nonumber}
\end{prop}

\subsection{Lower bound exploiting the bottleneck}
\label{subsec:low}

Consider the activity process $(X_t)_{t\geq 0}$ with activation rate $\nu$. 
For $S \subseteq \Omega$, let $\pi(S):=\sum_{(k,l) \in S } \pi_{(k,l)}(\nu)$ be the stationary probability of $S$. Define the {\it flow rate out of} $S$ as $Q(S,S^c):=\sum_{(k,l) \in S, (j,m) \in S^c} \pi_{(k,l)}(\nu) q((k,l),(j,m))$ and the {\it conductance of} $S$ as
$\Phi(S):=Q(S,S^c)/\pi(S).$

The conductance of the process $(X_t)_{t\geq 0}$ is defined as 
\[\Phi_*:=\min_{S \,: \, \pi(S) \leq 1/2} \Phi(S).\]

All the quantities we just defined depend on $\nu$, but we suppress it for conciseness.

The following result, valid for any Markov process on a finite state space $\Omega$ with conductance $\Phi_*$, shows how the conductance of the process yields a lower bound on the mixing time. It is a continuous-time version of Theorem 7.3 in \cite{LPW09} and since the proof is quite similar, it is omitted.
\begin{lem} \label{lem:cond}
For $ \epsilon \in \left( 0, \frac{1}{4} \right)$, $ t_{{\rm mix}}(\epsilon) \geq \left( \frac{1}{2} - 2 \epsilon \right) \frac{1}{\Phi_*}.
$
\end{lem}

We now exploit the fact that our activity process on the state space $\Omega$ has a geometric feature usually referred to as {\it bottleneck}, that strongly influences the mixing time. This bottleneck indeed makes some parts of $\Omega$ difficult to reach, resulting in a small conductance. As it turns out, $\C_2$ will be the bottleneck.

\begin{prop}
The conductance of $\C_2$ satisfies
\[
\Phi(\C_2) \sim L_2 \, \nu^{1-L_2}, \quad \ninf,
\]
and hence, for $\epsilon \in \left(0,\frac{1}{4}\right)$,
\eqn{ \label{lowbtmix}
t_{{\rm mix}}(\epsilon,\nu) \geq \Big( \frac{1}{2} - 2 \epsilon \Big) \frac{1}{L_2} \, \nu^{L_2-1}, \quad \ninf.
}
\end{prop}

\begin{IEEEproof}
Since $\C_1 \cup \{0\}$ will have at least half of the probability mass for $\nu$ sufficiently large, it is clear that $\pi(\C_2) \leq 1/2$ when $\ninf$. From \eqref{statdist} it follows that if $l$ and $m$ belong to the same component $k$ of size $L$, then
$
\frac{\pi_{(k,m)}(\nu)}{\pi_{(k,l)}(\nu)} = \frac{l! \,(L-l)!}{m! \, (L-m)!} \, \nu^{m-l},$ as $ \ninf.
$
Thus the conductance of $\C_2$ satisfies 
\[
\Phi(\C_2) = \frac{\pi_{(2,1)}(\nu) \cdot 1 }{ \sum_{l=1}^{L_2} \pi_{(2,l)}(\nu)}= \frac{ \frac{\pi_{(2,1)}(\nu)}{\pi_{(2,L_2)}(\nu)}}{ \sum_{l=1}^{L_2} \frac{\pi_{(2,l)}(\nu)}{\pi_{(2,L_2)}(\nu)}} \sim L_2 \, \nu^{1-L_2}.
\]
Then Lemma~\ref{lem:cond} gives the lower bound \eqref{lowbtmix},
since by definition $ \Phi_* \leq \Phi(\C_2)$.
\end{IEEEproof}

\section{Conclusion and extensions}
\label{sec6}

We have examined transient throughput characteristics
and associated starvation effects in CSMA networks in terms of the
transition times between dominant activity states.
We established how the magnitude of the transition time scales with
the activation rate and the sizes of the various network components
in partite interference graphs.
We also proved that in several cases the scaled transition time has an asymptotically exponential distribution and discussed the connection with
related exponentiality results for rare events and meta-stability
effects in statistical physics.
In addition, we investigated the convergence rate to equilibrium of the activity process in terms of mixing times.

In the present paper we have focused on partite interference graphs
with uniform activation rates, giving rise to a star-shaped state space,
but most of the methods and results extend to more general networks,
as long as there is a unique path between any two activity states.
For example, we could allow for the branches to be general trees,
with $M_{k,l}$ nodes in the $k$-th tree at distance~$l$ from the root,
and transition rates $b_{k,l} f_{k,l}(\nu)$ and $d_{k,l}$ away from
and towards the root, respectively.
The node-dependent functions $f_{k,l}(\nu)$ are particularly
relevant, since taking $f_{k,l}(\nu) = \nu^{1 / L_k}$ could for example
serve to balance the long-term throughputs in the various components.
In those cases, the expectation of the transition time will scale
differently, but most of the distributional results for the scaled
transition time continue to hold, even though qualitatively
different scenarios arise, with some branches being visited for
relatively short periods but an overwhelmingly large number of times.

\bibliographystyle{IEEEtran}
\bibliography{research} 

\appendix

\subsection{Proof of Lemma~\ref{asymptotic1}}
Order the state space as $\Omega=\{L,L-1,\dots, 1, 0\}$ and consider the generator matrix $Q(\nu)$ of the process $(X_t)_{t\geq 0}$ with $0$ an absorbing state. That is,
\[ Q(\nu)=
\begin{bmatrix}
 q_{L}(\nu) & L& 0 & \\
 \nu & q_{L-1}(\nu) & L-1 &  &   \\

 & \ddots & \ddots & \ddots &  \\

  &  & (L-1) \nu  &  q_{1}(\nu) & 1\\
    &  & 0 & 0 & 0\\
\end{bmatrix},
\]
where the diagonal elements are $q_l(\nu)=-(L-l) \nu -l$ for $l=1,\dots, L$. Write $Q(\nu)$ as
\[
Q(\nu)=\left(\begin{array}{cc}\mathbf{T}(\nu) & \mathbf{t}(\nu) \\ \mathbf{0} & 0 \\\end{array}\right), \nonumber
\]
where $\mathbf{T}(\nu)$ is an $L \times L$ invertible matrix. Since the characteristic polynomials of $-Q(\nu)$ and $-\mathbf{T}(\nu)$ satisfy the relation $ p_{-Q(\nu)}(z)= - z \, p_{-\mathbf{T}(\nu)}(z)$, the spectrum of $-Q(\nu)$ consists of that of $-\mathbf{T}(\nu)$ plus the eigenvalue zero with multiplicity one. Denote by $D(\nu)$ the $L\times L$ diagonal matrix, whose diagonal entries are $\{\sqrt{\theta_l(\nu)} \}_{i=L}^1$, where the $\theta$'s are the so-called \textit{potential coefficients}, defined as
$\theta_L(\nu)=1$ and $\theta_{l-1}(\nu)=\frac{l}{(L-l)\nu}\, \theta_{l}(\nu).$
The $L\times L$ matrix $G(\nu)=- D(\nu)^{1/2} \, \mathbf{T}(\nu) \, D(\nu)^{-1/2}$ is tridiagonal and symmetric with diagonal entries $g_{l,l}(\nu)=q_{L-l+1}(\nu)$ and $ g_{l,l+1}(\nu)=g_{l+1,l}(\nu)=-\sqrt{l \cdot (L-l+1) \nu}$. Since $G(\nu)$ is similar to $-\mathbf{T}(\nu)$, they have the same spectrum. Denote by $\mathcal{D}(p,R)$ the closed disc centered in $p$ with radius $R$, i.e.~$\mathcal{D}(p,R)=\{z \in \mathbb{C} ~:~ |z-p| \leq R\}.$
Consider the \textit{Gershgorin discs} $ \{\mathcal{D}_l(\nu)\}_{l=1}^L$ of $G(\nu)$, defined as $\mathcal{D}_l(\nu):=\mathcal{D}(-q_l(\nu),R_l(\nu))$, where the radius $R_l(\nu)$ is the sum of the absolute values of the non-diagonal entries in the $L-l+1$-th row, i.e.~$R_l(\nu):=\sum_{m \neq L-l+1} |g_{L-l+1,m}(\nu)|$.
Then
\eqan{
\mathcal{D}_{L}(\nu)&=\mathcal{D}(L,\sqrt{L \nu }),\nonumber\\
\mathcal{D}_{L-1}(\nu)&=\mathcal{D}(L-1+ \nu,\sqrt{L \nu}+\sqrt{ 2 (L-1) \nu}), \nonumber\\
& \, \, \, \, \dots  \nonumber\\
\mathcal{D}_{2}(\nu)&=\mathcal{D}(2+(L-2) \nu,\sqrt{ 3(L-2) \nu}+\sqrt{2(L-1) \nu}),\nonumber\\
\mathcal{D}_{1}(\nu)&=\mathcal{D}(1+(L-1) \nu,\sqrt{2(L-1) \nu}).\nonumber
}

We now exploit the second Gershgorin circle theorem, which is reproduced here for completeness.
\vspace{-0.47cm}\\
\begin{thma}
If the union of $j$ Gershgorin discs of a real $r\times r$ matrix $A$ is disjoint from the union of the other $r - j$ Gershgorin discs, then the former union contains exactly $j$ and the latter the remaining $r - j$ eigenvalues of $A$.\vspace{-0.47cm}\\
\end{thma}

In our case, for $\nu$ sufficiently large, the disc $\mathcal{D}_{L}(\nu)$ does not intersect with the union $\bigcup_{l=1}^{L-1} \mathcal{D}_l(\nu)$, thus the smallest eigenvalue $\alpha_1(\nu)$ lies in $ \mathcal{D}_{L}(\nu)$ and the other $L-1$ ones in $\bigcup_{l=1}^{L-1} \mathcal{D}_l(\nu)$. Hence, for $\nu$ sufficiently large, 
$ \alpha_1(\nu) \leq L+\sqrt{L \nu}$ and $ \alpha_i(\nu) \geq (L-1)+\nu-\sqrt{\nu}(\sqrt{L}+\sqrt{2(L-1)})$ for $i=2,\dots,L$. Therefore,
$ 0 < \frac{ \alpha_1(\nu)}{ \alpha_i(\nu)} \leq \frac{L+\sqrt{L \nu}}{\nu-\sqrt{\nu}(\sqrt{L}+\sqrt{2(L-1)})},$
and so $ \limnu \alpha_1(\nu) / \alpha_i(\nu) =0 $ for $i=2,\dots,L$. Hence,
\vspace{-0.29cm}\\
\[
\E T_{L,0}(\nu) \cdot \alpha_1(\nu)=1+\sum_{i=2}^{L} \frac{\alpha_1(\nu)}{ \alpha_i(\nu)} \to 1, \quad \ninf,
\]
\vspace{-0.3cm}\\
while for $ 2\leq  i \leq L$,
\vspace{-0.29cm}\\
\[
\E T_{L,0}(\nu) \cdot \alpha_i(\nu) > \frac{\alpha_i(\nu)}{\alpha_1(\nu)} \to \infty, \quad \ninf.
\]
\vspace{-0.3cm}

\subsection{Proof of Proposition~\ref{asyfinitemean}}

Using $\E S_M(\nu) = \E M \cdot \E X(\nu)$, gives
\[
\frac{S_M(\nu)}{\E S_M(\nu)} =
\frac{\sum_{i=1}^{M} X_i(\nu)}{\E M \E X(\nu)} =
\frac{1}{\E M} \sum_{i=1}^{M} \frac{X_i(\nu)}{\E X_i(\nu)}.
\]
Thus the Laplace transform of $S_M/ \E S_M$ may be written as
\eqan{
\LT_{S_M(\nu)/ \E S_M(\nu)}(s)
&= \LT_{\sum_{i=1}^{M} X_i(\nu) / \E X_i(\nu)}
\Big(\frac{s}{\E M}\Big) \nonumber \\
&= G_M \left(\LT_{X_i(\nu) / \E X_i(\nu)} \Big(\frac{s}{\E M}\Big)\right), \nonumber 
}
where $G_M(z)=\E (z^M)$. The assumption $X(\nu) / \E X(\nu) \cd Y$ implies that $\limnu \LT_{X_i(\nu) / \E X_i(\nu)}(t) = \LT_{Y}(t)$ for all $t \geq 0$ and $i \in \mathbb N$. Hence
\eqan{
\limnu \LT_{S_M(\nu)/ \E S_M(\nu)}(s)
&= \limnu G_M \left(\LT_{X_i(\nu) / \E X_i(\nu)}
\left(\frac{s}{\E M}\right)\right) \nonumber \\
&= G_M \left(\LT_{Y} \left(\frac{s}{\E M}\right)\right), \nonumber
}
which is the Laplace transform of $\frac{1}{\E M} \sum_{i=1}^{M} Y_i$. Invoking the continuity theorem for Laplace transforms concludes the proof.

\end{document}